\documentclass[10pt]{amsart}
\usepackage{amssymb}

\usepackage{amscd}
\usepackage{amsthm}
\usepackage[dvips]{graphicx}
\usepackage{color}
\usepackage[all]{xy}
\usepackage{verbatim}

\newtheorem{Lemma1}{{Lemma}}
\newtheorem{Theo1}[Lemma1]{{Theorem}}
\newtheorem*{Theo2}{{Theorem}}
\newtheorem{Def1}[Lemma1]{{Definition}}
\newtheorem{Prop1}[Lemma1]{{Proposition}}
\newtheorem{Claim1}[Lemma1]{{Claim}}
\newtheorem{Rem1}[Lemma1]{{Remark}}
\newtheorem{Cor1}[Lemma1]{{Corollary}}
\newtheorem{Ex1}[Lemma1]{{Example}}

\newenvironment{Lemma}{\begin{Lemma1}}{\end{Lemma1}}
\newenvironment{Def}{\begin{Def1}\em}{\end{Def1}}
\newenvironment{Prop}{\begin{Prop1}}{\end{Prop1}}
\newenvironment{Rem}{\begin{Rem1}\rm}{\end{Rem1}}
\newenvironment{Theorem}{\begin{Theo1}}{\end{Theo1}}

\newenvironment{Cor}{\begin{Cor1}}{\end{Cor1}}

\newenvironment{Example}{\begin{Ex1}\em}{\end{Ex1}}

\title{Degenerating $0$ in Triangulated Categories}

\author{Manuel Saor\'\i{}n}
\address{\newline
Departemento de Matem\'aticas,
\newline Universidad de Murcia, Aptdo. 4021
\newline 30100 Espinardo, Murcia,
\newline Spain}
\email{msaorinc@um.es}

\author{Alexander Zimmermann}
\address{\newline
Universit\'e de Picardie,
\newline D\'epartement de Math\'ematiques et LAMFA (UMR 7352 du CNRS),
\newline 33 rue St Leu,
\newline F-80039 Amiens Cedex 1,
\newline France}
\email{alexander.zimmermann@u-picardie.fr}

\date{January 28, 2019}

\newcommand{\uar}{\uparrow}

\newcommand{\lra}{\longrightarrow}
\newcommand{\lla}{\longleftarrow}
\newcommand{\ra}{\rightarrow}
\newcommand{\sdp}{\times\kern-.2em\vrule height1.1ex depth-.05ex}
\newcommand{\epi}{\lra \kern-.8em\ra}

\newcommand{\Z}{{\mathbb Z}}

\newcommand{\Hom}{\textup{Hom}}

 \normalbaselines
\begin{document}

\begin{abstract}
In previous work, based on work of Zwara and Yoshino, we defined and studied
degenerations of objects in triangulated categories analogous to
degeneration of modules. In triangulated categories it is surprising that the
zero object may degenerate. We study this systematically. In particular we show that
the degeneration of the zero object actually induces all other degenerations
by homotopy pullback, that degeneration of $0$ is closely linked, but not equivalent, to
having zero image in the Grothendieck group.
\end{abstract}

\maketitle

\section*{Introduction}
Degeneration of modules were intensively studied by e.g. Gabriel \cite{Gabriel},
Huisgen-Zimmermann, Riedtmann \cite{Riedtmann1}, Zwara \cite{Zwara,Zwara1}
since at least 1974, and was highly successful in various constructions. Degeneration of modules
is defined by the following setting.
Let $k$ be an algebraically closed field, and let $A$ be a finite dimensional $k$-algebra.
Then the $A$-module structures on the vector space $k^d$ form an affine algebraic variety
$\textup{mod}(A,d)$ on which
$GL_d(k)$ acts by conjugation. Isomorphism classes correspond to orbits under this action
and an $A$-module $M$ degenerates to $N$ if the point corresponding to $N$ belongs to the
Zariski closure of the $GL_d(k)$-orbit of the point corresponding to $M$. We write
$M\leq_{\textup{deg}}N$ in this case. Riedtmann and Zwara showed that $M\leq_{\textup{deg}}N$ if and
only if there is an $A$-module $Z$ and a short exact sequence
$0\ra Z\ra Z\oplus M\ra N\ra 0$. The last relation is denoted by $M\leq_{\textup{Zwara}} N$.
In collaboration with Jensen and Su \cite{JSZ1} the second named author started to study an analogous
concept for derived categories with a geometrically inspired concept based on orbit
closures, and then in \cite{JSZ2} more generally for triangulated categories based
on Zwara's characterisation replacing short exact sequences by distinguished triangles.
This last relation is denoted by the symbol $\leq_\Delta$.
Both concepts were highly successfully used in many places, cf e.g.
\cite{KellerYangZhou,KellerScherotzke,Krausecohom,Eisele1,Eisele2,Hiramatsu1,Hiramatsu2,Wang}.
Independently
Yoshino \cite{Yoshino} gave a scheme theoretic definition for degenerations in the
(triangulated) stable category of maximal Cohen-Macaulay modules, and he highlighted
that in $M\leq_\Delta N$ one should assume that the induced endomorphism on $Z$ should
be nilpotent. We denote the relation by $\leq_{\Delta+\textup{nil}}$ in this case.
Yoshino's scheme theoretic approach was a model for  us to give a more general geometric
definition for degeneration, which was achieved in \cite{SZ1} by introducing a scheme
theoretic degeneration $\leq_{\textup{cdeg}}$. We then showed that,
in case $\mathcal{T}$ has split idempotents,  $M\leq_{\textup{cdeg}}N$
always implies $M\leq_{\Delta+\text{nil}}N$, for objects $M,N\in\mathcal{T}$,
the converse being also true when $\mathcal{T}$ is the subcategory of compact
objects of a compactly generated algebraic triangulated category. We gave the
rather technical definition of $\leq_{\textup{cdeg}}$ by means of an
ambient triangulated category and a natural transformation $t$ of the identity
functor of this category satisfying a certain number of axioms. The actual
degeneration $M\leq_{\textup{cdeg}}N$ is then given by an object $Q$ in this ambient
category such that the cone of $t_Q$ is $N$, and such that $M$ and $Q$ become
isomorphic when one inverts all morphisms $t_X$  for all
objects $X$ in the ambient category
in the Gabriel-Zisman sense.
For more ample details we refer to Section~\ref{veriousconceptsofdegen}.

A striking phenomenon is that, unlike in the module case, in triangulated
categories $\mathcal T$
one may have non zero objects $M$ with $0\leq_{\Delta+\textup{nil}}M$.
In the present paper we study this phenomenon systematically.

As a first result we show that $0\leq_{\textup{cdeg}} M$ is given by objects $Q$ so
that the identity transformation $t$ of the degeneration data has the property
that $t_Q$ is nilpotent.

Second, we show in the present paper that if $M\leq_{\Delta}N$ via a
distinguished triangle
$$\xymatrix{Z\ar[r]^-{\alpha\choose\beta}&Z\oplus M\ar[r]^-{(\gamma,\delta)}& N\ar[r]& Z[1],}$$
then $0\leq_\Delta \text{cone}(\alpha)$ and $0\leq_\Delta\textup{cone}(\delta)$, and if
$0\leq_{\Delta} X$ via a distinguished triangle
$$\xymatrix{Z\ar[r]^-\alpha& Z\ar[r]& X\ar[r]& Z[1],}$$
then for any morphism $M\lra Z$ the
homotopy pullback $N$
along $M\lra Z\stackrel{\alpha}\lla Z$ gives a degeneration $M\leq_{\Delta}N$.
Analogous statements hold for $\leq_{\Delta+\text{nil}}$.
Note that $0\leq_{\Delta}M$ if and only if $M$ is the cone of an endomorphism $v$
of some object $Z$, and $0\leq_{\Delta+\textup{nil}}M$ if and only if $v$ is in
addition nilpotent.  Hence, degeneration of $0$ is very intrinsic in
degenerations of triangulated categories. It is not possible to get rid of
this phenomenon, when one wants to study degenerations.  Therefore, our
observation interprets degeneration as some sort of deformation along a (nilpotent)
endomorphism of some parameter space $Z$ via a homotopy pullback construction.

Third, it is quite clear that $0\leq_\Delta N$ implies that $N$ has vanishing
image in the Grothendieck group of the triangulated category. We characterise
the objects with image being the zero object in the Grothendieck group by showing
that this condition is equivalent to the fact that
$M\leq_\Delta \bigoplus_{i=1}^r(X_i\oplus X_i[t_i])$ for pairwise different  odd integers
$t_i$, and objects $X_i$ of $\mathcal T$. This indicates already that the objects
$X\oplus X[t]$, with $t$ odd, are very intimately linked to degeneration phenomena.
When $X$ runs through a generating set of the triangulated category, we show that
the triangulated category generated by these objects coincides with the full
triangulated subcategory $\mathcal T^0$ of $\mathcal T$ consisting of objects with
image $0$ in the Grothendieck group. We further show that the triangulated subcategory
generated by objects which are degenerations of $0$ coincides with ${\mathcal T}^0$.
We give a class of examples of a triangulated category $\mathcal T$,
and an object $X$ in $\mathcal T$ such that
$[X]=0$ in $K_0({\mathcal T})$, but $0\not\preceq_{\Delta+\textup{ nil}}X$, where
$\preceq_{\Delta+\textup{ nil}}$
is the transitive hull of the relation $\leq_{\Delta+\textup{nil}}.$
It is not hard to get another class of examples where $[X]=0$ in $K_0({\mathcal T})$
but where $0\not\leq_{\Delta}X$.

The paper is organised as follows. In Section~\ref{veriousconceptsofdegen}
we recall the necessary concepts on the various types of degeneration.
In Section~\ref{cdegeneration} we study $\leq_{\textup{cdeg}}$
and show that $0\leq_{\textup{cdeg}}M$ if and only if $t_Q$ is nilpotent for the
corresponding degeneration data.
In Section~\ref{degenerationofzeroimpliesalltheothers} we show that triangle
degeneration of zero is always present in triangle degenerations in
triangulated categories. Indeed all triangle degenerations are obtained from degenerations of zero via homotopy cartesian squares. In Section~\ref{ZeroinGrothendieck} we study the image of
triangle degenerations of $0$ in the Grothendieck group.

\section{Review on Degenerations in Triangulated Categories}

\label{veriousconceptsofdegen}

We have different degeneration concepts. The first one, the triangle degeneration, is a triangular
category analogue of Zwara's definition of degeneration in the case of module categories. Zwara
says \cite{Zwara,Zwara1}
that for a $k$-algebra $A$ an $A$-module $M$ degenerates to an $A$-module $N$ if
and only if there is an $A$-module $Z$ and a short exact sequence $0\ra Z\ra Z\oplus M\ra N\ra 0$.
Yoshino~\cite{Yoshino} highlighted the importance of assuming that the induced endomorphism of
$Z$ should be nilpotent. In case of a category where Fitting's lemma holds we can always assume this
fact.

\begin{Def} \cite{JSZ2,SZ1}\label{deltadegen}
Let $K$ be a commutative ring and let $\mathcal T$ be a $K$-linear triangulated category.
Then for two objects $M$ and $N$ in $\mathcal T$ we get $M\leq_{\Delta}N$, accasionally written  $M\leq_{\Delta,\textup{left}}N$,  if and
only if there is $Z$ and a distinguished triangle
$$\xymatrix{Z\ar[r]^-{v\choose u}&Z\oplus M\ar[r]& N\ar[r]& Z[1].}$$
We say that $M\leq_{\Delta+\textup{nil}}N$, or occasionally
$M\leq_{\Delta+\textup{nil, left}}N$, if and only if
there is such a distinguished triangle with $v$ is nilpotent.
\end{Def}

Note that by \cite{SZ2} $M\leq_{\Delta+\textup{nil}}N$ implies that there is an
object $Z'$ and a distinguished triangle
$$\xymatrix{N\ar[r]& Z'\oplus M\ar[r]^-{(v',u')}& Z'\ar[r]& N[1].}$$

We write   $M\leq_{\Delta ,\textup{right}}N$ (resp. $M\leq_{\Delta+\textup{nil, right}}N$) if there is such a distinguished triangle
(with $v'$ nilpotent). Note that $M\leq_{\Delta,\textup{left}}N$ (resp. $M\leq_{\Delta+\textup{nil,left}}N$) in $\mathcal{T}$ if and only if  $M\leq_{\Delta,\textup{right}}N$ (resp. $M\leq_{\Delta+\textup{nil,right}}N$)  in the opposite category $\mathcal{T}^{op}$. So categorical duality applies and results about $\leq_{\Delta,\textup{left}}$ (resp. $\leq_{\Delta+\textup{nil, left}}$) admit categorical dual ones, that we sometimes omit to state.

If $\mathcal T$ has
split idempotents and artinian endomorphism rings of objects, or if $\mathcal T$
is the category of compact objects in a compactly generated algebraic triangulated category, then
$M\leq_{\Delta+\textup{nil, right}}N$ if and only if
$M\leq_{\Delta+\textup{nil, left}}N$ (see \cite[Theorem 1]{SZ2}). For a general triangulated category $\mathcal{T}$, when we write $\leq_\Delta$ (resp. $\leq_{\Delta+\textup{nil}}$) we will mean  $\leq_{\Delta,\textup{left}}$ (resp. $\leq_{\Delta+\textup{nil, left}}$).

\begin{Cor}\label{K0triangledeg}
Let $K$ be a commutative ring and let $\mathcal T$ be a $K$-linear triangulated category.
Then $M\leq_{\Delta}N$ implies that $[M]=[N]$, where $[X]$ denotes the image of $X$ in $K_0({\mathcal T})$.
\end{Cor}

Indeed, this is a direct consequence of the fact that there is a distinguished triangle
$$\xymatrix{Z\ar[r]^-{v\choose u}&Z\oplus M\ar[r]& N\ar[r]& Z[1].}$$
The distinguished triangle shows $[Z]+[N]-[Z\oplus M]=0$ in the Grothendieck group, and therefore
$[M]=[N]$.

\medskip

A second concept of degeneration, motivated by Yoshino's work, is given by the following definition.

\begin{Def}\cite{SZ1} \label{degendatadef}
Let $K$ be a commutative ring and let ${\mathcal C}_K^\circ$ be a
$K$-linear triangulated category with split idempotents.

A degeneration data for ${\mathcal C}_K^\circ$ is given by
\begin{itemize}
\item a triangulated category ${\mathcal C}_K$ with split idempotents and a fully faithful
embedding ${\mathcal C}_K^\circ\longrightarrow{\mathcal C}_K$,
\item
a triangulated category   ${\mathcal C}_V$ with split idempotents
and a full triangulated subcategory ${\mathcal C}_V^\circ$,
\item  triangulated functors $\uar_K^V:{\mathcal C}_K\lra {\mathcal C}_V$,
which we write after the arguments,
and  $\Phi:{\mathcal C}_V^\circ\ra {\mathcal C}_K$, so that
$({\mathcal C}_K^\circ)\uar_K^V\subseteq {\mathcal C}_V^\circ$, when
we view ${\mathcal C}_K^\circ$ as a full subcategory of ${\mathcal
C}_K$,
\item  a natural transformation
$\text{id}_{{\mathcal C}_V}\stackrel{t}{\lra} \text{id}_{{\mathcal C}_V}$
of triangulated functors such that
\item
for each object $M$ of ${\mathcal C}_K^\circ$ the morphism
$\xymatrix{\Phi(M\uar_K^V)\ar[rr]^-{\Phi(t_{M\uar_K^V})}&&\Phi(M\uar_K^V)}$
is a split monomorphism in ${\mathcal C}_K$ with cone $M$.
\end{itemize}
\end{Def}

Degeneration is then given by the following concept.

\begin{Def}\cite{SZ1} \label{degendef}
Given two objects $M$ and $N$ of ${\mathcal C}_K^\circ$ we say that
$M$ degenerates to $N$ in the categorical sense if there is a
degeneration data for ${\mathcal C}_K^\circ$ and an object $Q$ of
${\mathcal C}_V^\circ$ such that
$$p(Q)\simeq p(M\uar_K^V)\mbox{ in ${\mathcal C}_V^\circ[t^{-1}]$
and }\Phi(\textup{cone}(t_Q))\simeq N,$$
 where ${\mathcal C}^\circ_V[t^{-1}]$ is
the Gabriel-Zisman localisation at the endomorphisms $t_X$ for all
objects $X$ of  ${\mathcal C}^\circ_V$, and where $p:{\mathcal
C}_V^\circ\longrightarrow{\mathcal C}_V^\circ [t^{-1}]$ is the
canonical functor. In this case we write $M\leq_{\textup{cdeg}}N$.
\end{Def}

\begin{Rem}
We note that this concept is a generalisation to general triangulated categories
of a definition given by Yoshino~\cite{Yoshino} for the case of
stable categories of maximal Cohen-Macaulay
modules over a local Gorenstein algebra.
\end{Rem}

We get the following connection on these two concepts
$\leq_{\Delta+\textup{nil}}$ and $\leq_{\textup{cdeg}}$.

\begin{Theorem}\cite{SZ1}\label{SZ1thm}
Let $\mathcal T$ be a $K$-linear triangulated category with split idempotents.
Then $$M\leq_{\textup{cdeg}}N\Rightarrow M\leq_{\Delta+\textup{nil}}N$$ and
if $\mathcal T$
is the category of compact objects in a compactly generated algebraic triangulated category
then $$M\leq_{\textup{cdeg}}N\Leftrightarrow M\leq_{\Delta+\textup{nil}}N.$$
\end{Theorem}

\begin{Cor}\label{K0catdeg}
Let $\mathcal T$ be a $K$-linear triangulated category with split idempotents.
Then $$M\leq_{\textup{cdeg}}N\Rightarrow [M]=[N]$$ in $K_0({\mathcal T})$.
\end{Cor}

Indeed, this is Theorem~\ref{SZ1thm} in connection with Corollary~\ref{K0triangledeg}.

\medskip

Let $\mathcal T$ be a skeletally small triangulated category.
In \cite{JSZ2} it is studied when $\leq_\Delta$ is a partial order on the isomorphism
classes of objects in $\mathcal T$, and in particular when
$\leq_\Delta$ is a transitive relation. In general this is not the case.
However  the following was proved there.

\begin{Prop} \cite[Proposition 2]{JSZ2} \label{JSZ2prop}
Let $\mathcal T$ be a triangulated category with split idempotents such that the
endomorphism ring of each object in $\mathcal T$ is artinian.
Then $\leq_{\Delta+\textup{ nil}}$ coincides with $\leq_{\Delta}$ and both relations are
transitive and reflexive on the set of isomorphism classes of objects in $\mathcal T$.
\end{Prop}

\section{Scheme theoretic degeneration of zero}

\label{cdegeneration}

\subsection{Torsion Degeneration Data}

Recall that in \cite{Yoshinomodules} Yoshino defines as well the scheme theoretic
definition of a degeneration of modules. Let $A$ be a $k$-algebra for
some field $k$.
Yoshino says that the $A$-module $M$ degenerates to the $A$-module $N$
along a discrete valuation ring if
there is a discrete valuation $k$-algebra $V$ with uniformiser $t$ and a finitely
generated $A\otimes_kV$-module $Q$ which is flat as $V$-module, such that
$Q/tQ\simeq N$, and such that $Q[t^{-1}]\simeq M\otimes_kV[t^{-1}]$. In
\cite{Yoshino} a similar setting was used for degeneration in a triangulated
(stable) category, and there the condition on $Q$ to be flat over $V$ is missing. We now
explain why this is the case.

\begin{Prop}\label{geometricdegenerationfortriangcat}
Let ${\mathcal C}_K^\circ$ be a $K$-linear triangulated category with split idempotents.
Then $0\leq_{\textup{cdeg}}N$ if and only if there is a degeneration data
$({\mathcal C}_K^\circ,{\mathcal C}_K,{\mathcal C}_V^\circ,{\mathcal C}_V,\uar_K^V,\Phi,t)$
and $Q$ an object
in ${\mathcal C}_V^\circ$ such that $t_{Q}$ is nilpotent with
$$\textup{cone}(\Phi(t_Q))\simeq \Phi(\textup{cone}(t_Q))\simeq N\mbox{ and }\Phi(Q)\in {\mathcal C}_K^\circ.$$
\end{Prop}

\begin{proof}
Suppose we have a degeneration data and by \cite[Lemma 7]{SZ1}
$Q$ an object in ${\mathcal C}_V^\circ$
so that $t_{Q}$ is nilpotent with
$$\textup{cone}(\Phi(t_Q))\simeq \Phi(\textup{cone}(t_Q))=:N\in {\mathcal C}_K^\circ.$$
Then $p(Q)\simeq 0$ since $t_Q$ is nilpotent. Hence, by definition
$0\leq_{\textup{cdeg}}N$.

If $0\leq_{\textup{cdeg}}N$, then there is a degeneration data and an object $Q$
in ${\mathcal C}_V^\circ$ with nilpotent $t_Q$ and
$\Phi(\textup{cone}(t_Q))=N\in {\mathcal C}_K^\circ.$ Note that the proof of
the first part of Theorem~\ref{SZ1thm} gives that actually $\Phi(\textup{cone}(f))$
is the object named $Z$, where
$f:M\uar_K^V\lra Q$ is a morphism such that $p(f)$ is an isomorphism. In case $M=0$,
then $\textup{cone}(f)=Q$, and hence $\Phi(Q)=Z\in\mathcal{C}_K^\circ.$
\end{proof}

\begin{Rem}
We just proved that when we want to consider degenerations $M\leq_{\textup{cdeg}} N$
without direct factor $0\leq_{\textup{cdeg}}X$, then it is necessary and sufficient to
add as an additional condition that $Q$ does not have any non zero direct
factor $Q'$ with $\Phi(Q')\in\mathcal{C}_K^\circ$ such that $t_{Q'}$ is nilpotent.
\end{Rem}

\subsection{Triangle degeneration of zero}

\begin{Lemma}
Let $\mathcal T$ be a triangulated category with split idempotents. Then
$0\leq_{\Delta}N$ (resp.
$0\leq_{\Delta+\textup{nil}}N$) is equivalent to
$N=\textup{cone}(v)$ for some (resp. nilpotent) endomorphism $v$
of an object $Z$. Moreover, in this case
$N\leq_{\Delta+\textup{nil}} Z\oplus Z[1]$.
\end{Lemma}

\begin{proof}
If $\mathcal T$ is a $K$-linear triangulated category with split idempotents.
Then $0\leq_{\Delta}N$  (resp.
$0\leq_{\Delta+\textup{nil}}N$) is equivalent to the existence of an object $Z$ of $\mathcal T$
and a  (resp. nilpotent) endomorphism $v$ of $Z$ such that $N=\textup{cone}(v)$. In
other words there is a distinguished triangle
$$\xymatrix{Z\ar[r]^-v& Z\ar[r]^-{\epsilon}& N\ar[r]^-{\delta}& Z[1].}$$
Then, since $$\xymatrix{0\ar[r]& Z\ar[r]^-{\text{id}_Z}&Z\ar[r]& 0}$$
is a distinguished triangle, also
$$
\xymatrix{Z\oplus Z[-1]\ar[rr]^-{\scriptsize\left(\begin{array}{cc}v&0\end{array}\right)}&&
Z\ar[r]^-{\scriptsize\left(\begin{array}{c}\epsilon\\ 0\end{array}\right)}&
N\oplus Z\ar[rr]^-{\scriptsize\left(\begin{array}{cc}\delta&0\\0& \text{id}_Z\end{array}\right)}&&Z[1]\oplus Z}
$$
is a distinguished triangle. Shift to the left gives a distinguished triangle
$$
\xymatrix{
Z\ar[r]^-{\scriptsize\left(\begin{array}{c}\epsilon\\ 0\end{array}\right)}&
N\oplus Z\ar[rr]^-{\scriptsize\left(\begin{array}{cc}\delta&0\\0& \text{id}_Z\end{array}\right)}&&Z[1]\oplus Z
\ar[rr]^-{\scriptsize\left(\begin{array}{cc}v&0\end{array}\right)}&&Z[1]}
$$
Now, $0$ is clearly nilpotent, and hence $N\leq_{\Delta+\text{nil}}Z\oplus Z[1]$).
\end{proof}

\begin{Theorem} \label{degenerationshomotopycartesiansquare}
Let $\mathcal T$ be a triangulated $k$-category. The following assertions hold:
\begin{enumerate}
\item If $M\leq_{\Delta,\textup{left}}N$ (resp. $M\leq_{\Delta+\textup{nil,left}}N$) via a distinguished triangle $$\xymatrix{
Z\ar[r]^-{\scriptsize\left(\begin{array}cu\\v\end{array}\right)}&M\oplus Z\ar[rr]^-{\scriptsize\left(\begin{array}{cc}\pi&\tau\end{array}\right)}&&N\ar[r]^{\mu}&Z[1]
},$$ where $v$ is a nilpotent endomorphism for $\leq_{\Delta+\textup{nil,left}}$, then $0\leq_{\Delta,\textup{left}}\textup{cone}(\pi)\simeq\textup{cone}(v)$ (resp. $0\leq_{\Delta+\textup{nil, left}}\textup{cone}(\pi)\simeq\textup{cone}(v)$).
\item If $M\leq_{\Delta,\textup{right}}N$ (resp. $M\leq_{\Delta+\textup{nil,right}}N$) via a distinguished triangle $$\xymatrix{
N\ar[r]^-{\scriptsize\left(\begin{array}c\sigma\\ \lambda\end{array}\right)}&M\oplus Z'\ar[rr]^-{\scriptsize\left(\begin{array}{cc} u'& v'\end{array}\right)}&&Z'\ar[r]^{\mu'}&N[1]
},$$ where $v'$ is a nilpotent endomorphism for $\leq_{\Delta+\textup{nil,right}}$, then $0\leq_{\Delta,\textup{right}}\textup{cone}(\sigma)\simeq\textup{cone}(v')$ (resp. $0\leq_{\Delta+\textup{nil, right}}\textup{cone}(\sigma)\simeq\textup{cone}(v')$).
\end{enumerate}

On the other hand,  if $C$ is the cone of the (resp. nilpotent) endomorphism $v$ of $Z$, then
for every morphism $u:M\lra Z$, by homotopy pullback along $u$ we obtain $N$
such that  $M\leq_{\Delta, \textup{right}}N$ (resp. $M\leq_{\Delta+\textup{nil, right}}N$ when $v$ is nilpotent).  Similarly, for every $u':Z\lra M$ by homotopy pushout
along $u'$ we obtain $N$ such that  $M\leq_{\Delta, \textup{left}}N$ (resp. $M\leq_{\Delta+\textup{nil, left}}N$ when $v'$ is nilpotent).
\end{Theorem}
\begin{proof}
Assertion 2 is the categorical  dual of assertion 1. So we just prove this last assertion. Consider the distinguished triangle $$\xymatrix{
Z\ar[r]^-{\scriptsize\left(\begin{array}{c}u\\ v\end{array}\right)}&M\oplus Z\ar[rr]^{\scriptsize\left(\begin{array}{cc}
\pi&\tau\end{array}\right)}&&N\ar[r]&Z[1]
}. $$ By definition, the square $$
\xymatrix{
Z\ar[r]^v\ar[d]^u&Z\ar[d]^\tau\\
M\ar[r]^\pi&N
}
$$ is homotopy cartesian. We then complete the horizontal maps to distinguished triangles to get a morphism between triangles.
$$
\xymatrix{
Z\ar[r]^v\ar[d]^u&Z\ar[d]^\tau\ar[r]^\nu&C(v)\ar[r]\ar@{-->}[d]^{\alpha}&Z[1]\ar[d]^{u[1]}\\
M\ar[r]^\pi&N\ar[r]^\mu&C(\pi)\ar[r]&M[1]
}.
$$  We may choose $\alpha$ to be an isomorphism and, by the upper row of the diagram, we conclude that  $0\leq_{\Delta,\textup{left}}C(v)\simeq C(\pi)$, resp.  $0\leq_{\Delta+\textup{nil, left}}C(v)\simeq C(\pi)$ when $v$ is nilpotent.

For the last statement, we just prove one half of it since the other half is dual.
Let us consider  the distinguished triangle
$$\xymatrix{Z\ar[r]^{v}&Z\ar[r]&C\ar[r]&Z[1]}.$$

For any morphism $u:M\lra Z$ we may complete the diagram
$$
\xymatrix{
&M\ar[d]^{u}&\\
Z\ar[r]^{v}&Z\ar[r]^\mu&C\ar[r]&Z[1]
}
$$
to
$$
\xymatrix{
&M\ar[d]^{u}\ar[r]^{\mu u}&C'\ar[d]^{\textup{id}_{C}}\\
Z\ar[r]^{v}&Z\ar[r]^\mu&C\ar[r]&Z[1]
}
$$
By the dual of Neeman~\cite[Lemma 1.4.3]{Neeman} there is an object $N$ and a morphism $v$ such that
$$
\xymatrix{N\ar[r]^{\sigma}\ar[d]^{v}
&M\ar[d]^{u}\ar[r]^{\mu u}&C\ar[d]^{\textup{id}_{C}}\ar[r]&N[1]\ar[d]^{v[1]}\\
Z\ar[r]^{v}&Z\ar[r]^\mu&C\ar[r]&Z[1]
}
$$
is a morphism of distinguished triangles such that
$$
\xymatrix{N\ar[r]^{\sigma}\ar[d]^{v}
&M\ar[d]^{u}\\
Z\ar[r]^{v}&Z
}
$$
is a homotopy cartesian square.
By definition this gives a distinguished triangle
$$
\xymatrix{N\ar[r]&M\oplus Z\ar[r]&Z\ar[r]&N[1]}
$$
Hence  $M\leq_{\Delta,\textup{right}}N$, (resp. $M\leq_{\Delta+\textup{nil, right}}N$ when $v$ is nilpotent).
\end{proof}

\begin{Def}
Let $\mathcal T$ be a triangulated category, let
$M$ and $Z$ be objects in $\mathcal T$ and let $v$ be an endomorphism of $Z$.
For $u\in \Hom_{\mathcal T}(M,Z)$ we denote by $(Deg(u,v),r,s)$,
or just $Deg(u,v)$ for
short, the homotopy pullback
$$\xymatrix{Z\ar[r]^v&Z\\ Deg(u,v)\ar[u]^r\ar[r]^-s&M\ar[u]^u.}$$
Dually, for $u'\in Hom_{\mathcal T}(Z,M)$ we denote by $(Ged(u',v),r,s)$,
or just $Ged(u',v)$ for
short, the homotopy pushout
$$\xymatrix{Z\ar[d]^{u'}\ar[r]^v&Z\ar[d]^r\\ M\ar[r]^-s&Ged(u',v).}$$
\end{Def}

We then have
$$M\leq_{\Delta,\textup{ right}}Deg(u,v)\text{ and }M\leq_{\Delta,\textup{ left}}Ged(u',v)$$
and if $v$ is in addition nilpotent, then
$$M\leq_{\Delta+\textup{nil,\textup{ right}}}Deg(u,v)\text{ and }M\leq_{\Delta+\textup{nil,\textup{ left}}}Ged(u',v).$$
 %

\begin{Rem}\label{SZ2intermsofGEDEG}
Suppose $\mathcal T$ has split idempotents and
either endomorphism rings of objects in $\mathcal T$ are artinian or else
$\mathcal T$ is the category of compact objects in an algebraic compactly
generated triangulated category. Then \cite{SZ2} shows that
\begin{itemize}
\item
for each object $Z$, nilpotent endomorphism $v$ of $Z$ and morphism $M\stackrel{u}\lra Z$
there is an object $Z'$, a nilpotent endomorphism $v'$ of $Z'$ and a morphism $Z'\stackrel{u'}\lra M$
with $Ged(u',v')\simeq Deg(u,v)$,
\item
and
for each object $Z'$, nilpotent endomorphism $v'$ of $Z'$ and morphism $Z'\stackrel{u'}\lra M$
there is an object $Z$, a nilpotent endomorphism $v$ of $Z$ and a morphism $M\stackrel{u}\lra Z$
with $Ged(u',v')\simeq Deg(u,v)$.
\end{itemize}
\end{Rem}

\begin{Example}
Consider a finite dimensional $k$-algebra $A$ over an algebraically closed field $k$.
Let $M$ and $N$ be two finite dimensional $A$-modules, and suppose
$M\leq_{\textup{deg}} N$.  Then by the Zwara-Riedtmann theorem \cite{Zwara1} there is a
finite dimensional $A$-module $Z$ and a short exact sequence
$$\xymatrix{0\ar[r]& Z\ar[r]^-{v\choose u}& Z\oplus M\ar[r]^-{(\tau,\pi)}& N\ar[r]& 0.}$$
A trivial case is when there is a short exact sequence
$$\xymatrix{0\ar[r]& N_1\ar[r]^-{\iota}& M\ar[r]^-{\rho}&N_2\ar[r]& 0.}$$
Taking $Z=N_1$ the short exact sequence
$$
\xymatrix{0\ar[r]&N_1\ar[r]^-{\scriptsize\left(\begin{array}{c}\iota\\0\end{array}\right)}&M\oplus N_1\ar[rr]^-{\scriptsize\left(\begin{array}{cc}\rho&0\\0&\textup{id}_{N_1}\end{array}\right)}&&N_2\oplus N_1\ar[r]&0}
$$
is a Zwara-Riedtmann sequence. In particular, in this case $v=0$.
Since $A-\textup{mod}\hookrightarrow \mathcal{D}^b(A)$, by mapping a module to the corresponding stalk
complex in degree $0$, we can consider the Zwara exact sequence as three terms
of a distinguished triangle. By Theorem~\ref{degenerationshomotopycartesiansquare}
the map $\pi$ in the general setting yields a
degeneration $0\leq_{\Delta+\textup{nil}}\textup{cone}(\pi)$ in $\mathcal{D}^b(A)$.
Of course there is no degeneration of $0$ of modules. Hence, this phenomenon is a purely
triangulated one. What does this give in  this particular situation?
Actually, for a general $\pi$, we get
$$\textup{cone}(\pi)=(\dots\lra 0\lra M\stackrel{\pi}\lra N\lra 0\lra\dots)$$
is the two term complex concentrated in degrees $-1$ and $0$. In the special
situation of the short exact sequence
$$\xymatrix{0\ar[r]& N_1\ar[r]^-{\iota}& M\ar[r]^-{\rho}&N_2\ar[r]& 0}$$
we see that $v=0$ and
$$\textup{cone}(\pi)=\textup{cone}\left(\begin{array}{cc}\rho&0\end{array}\right)
\simeq (\dots\ra M\stackrel{\rho}\ra N_2\ra 0\dots)\oplus N_1\simeq
\ker(\rho)[1]\oplus N_1\simeq N_1[1]\oplus N_1$$
since $\rho$ is surjective, and hence the two term complex given by $\rho$
is isomorphic, in the derived category, to the kernel of $\rho$.
This is coherent with the computation of the degeneration given by the zero map on $Z=N_1$.
\end{Example}

\section{Triangle degeneration as homotopy cartesian square: some consequences}

\label{degenerationofzeroimpliesalltheothers}

\subsection{The case of a single degeneration object}

In view of Theorem~\ref{degenerationshomotopycartesiansquare} we consider
degeneration as homotopy cartesian squares.

\begin{Lemma}\label{degenerationfactorises}
Let $\mathcal{T}$ be a triangulated category, let $Z\stackrel{w}{\longrightarrow}M\stackrel{u}{\longrightarrow}Z$ be morphisms in $\mathcal{T}$ and let $\nu_1,\nu_2:Z\longrightarrow Z$ be endomorphisms. Denote
$Deg(u,\nu_1)=(Deg(u,\nu_1),r_1,s_1)$ and $Ged(w,\nu_2)=(Ged(w,\nu_2),r_2,s_2)$ for short.
The following assertions hold:
\begin{enumerate}
\item $Deg(u,\nu_1)\leq_{\Delta, \textup{right}} Deg(u,\nu_1\nu_2)$ as well as $Deg(u,\nu_1\nu_2)=Deg(r_1,\nu_2)$.
When $\nu_2$ is nilpotent, we can replace $\leq_{\Delta, \textup{right}}$ by $\leq_{\Delta+\textup{nil,right}}$.
\item $Ged(w,\nu_2)\leq_{\Delta,\textup{left}} Ged(w,\nu_1\nu_2)$  as well as $Ged(u,\nu_1\nu_2)=Ged(r_2,\nu_1)$.
    When $\nu_1$ is nilpotent, we can replace $\leq_{\Delta,\textup{left}}$ by $\leq_{\Delta+\textup{nil,left}}$.
\end{enumerate}
\end{Lemma}

\begin{proof}
Assertion 1 follows from assertion 2 by categorical duality, so we just prove assertion 2. Let us consider the following   diagram, where the upper inner and the outer squares are homotopy cartesian
$$
\xymatrix{
Z\ar[r]^w\ar[d]_{\nu_2}&M\ar[d]^{s_2}\ar@/^3pc/[dd]^{s_{21}}\\
Z\ar[r]^-{r_2}\ar[d]_{\nu_1}&Ged(w,\nu_2)\ar@{-->}[d]^{f}\\
Z\ar[r]^-{r_{21}}&Ged(w,\nu_1\nu_2)
}
$$
There is then a dotted morphism $f:Ged(w,\nu_2)\longrightarrow Ged(w,\nu_1\nu_2)$ completing commutatively the diagram (see \cite[page 54]{Neeman}). But, by \cite[Lemma 9]{SZ2}, we can choose $f$ so that the square
$$
\xymatrix{
Z\ar[r]^-{r_2}\ar[d]_{\nu_1}&Ged(w,\nu_2)\ar[d]^f\\
Z\ar[r]^-{r_{21}}&Ged(w,\nu_1\nu_2)
}
$$
is homotopy cartesian. Therefore we have $Ged(w,\nu_2)\leq_{\Delta,\textup{left}} Ged(w,\nu_1\nu_2)$
as well as $Ged(u,\nu_1\nu_2)=Ged(r_2,\nu_1)$. We can replace $\leq_{\Delta,\textup{left}}$ by $\leq_{\Delta+\textup{nil,left}}$ when $\nu_1$ is nilpotent.
\end{proof}

%

\subsection{The case of two degeneration objects}

In this subsection we show that a degeneration $M\leq_{\Delta,\textup{right}}N$ obtained by two different triangles, but with the same morphism $N\lra M$, naturally yields another degeneration. Concretely:

\begin{Prop} \label{pushouts}
Let $\mathcal T$ be a triangulated category and let $M\leq_{\Delta,\textup{right}}N$ (resp. $M\leq_{\Delta+\textup{nil,right}}N$) be a degeneration obtained by the following homotopy cartesian squares

$$\xymatrix{
Z\ar[r]^v&Z\\
N\ar[u]^t\ar[r]^s&M\ar[u]^u}
\text{\;\;\;\; and \;\;\;\;}
\xymatrix{
Z'\ar[r]^{v'}&Z'\\
N\ar[u]^{t'}\ar[r]^{s'}&M\ar[u]^{u'}
}$$ where $v$ and $v'$ are assumed to be nilpotent endomorphism in the $\leq_{\Delta+\textup{nil,right}}$-case,
 and suppose that  $s=s'$.
If  $X$ denotes the (lower right corner of) the homotopy pushout along $Z\stackrel{t}{\lla} N\stackrel{t'}{\lra} Z'$ and $Y$ denotes that of the
 homotopy pushout along $Z\stackrel{u}{\lla} M\stackrel{u'}{\lra} Z'$, then $X\leq_{\Delta, \textup{left}}Y$ (resp.  $X\leq_{\Delta+ \textup{nil,left}}Y$) .
\end{Prop}

\begin{Rem}
Note that the order $X\leq_\Delta Y$ is inverse to the order $M\leq_\Delta N$.
\end{Rem}

\begin{Rem}
Note that Proposition~\ref{pushouts} should be seen in the context of
Bongartz \cite[Lemma 1.1]{Bongartz}. There it is shown that if
$$0\ra M'\ra M\ra M''\ra 0$$ is a short exact sequence, then $M$ degenerates to the
pushout $N_\lrcorner$ along an endomorphism of $M'$, and $M$ degenerates to the pullback
$N_\ulcorner$
along an endomorphism of $M''$.
\end{Rem}

\begin{proof} (of Proposition~\ref{pushouts}).
Note that $X$ and $Y$ are uniquely determined up to non-unique isomorphism. However,
we have some freedom for the choice of the morphisms which complete
$Z\stackrel{t}{\lla} N\stackrel{t'}{\lra} Z'$ and  $Z\stackrel{u}{\lla} M\stackrel{u'}{\lra} Z'$
to the corresponding homotopy pushouts.  We claim that we can choose the morphisms
$w:Z\longrightarrow Y$ and $w':Z'\longrightarrow Y$ so that the two squares and
the patching of them in the following commutative diagram are homotopy cartesian:

$$\xymatrix{N \ar[d]_{t} \ar[r]^{s} & M \ar[d]_{u} \ar[r]^{u'} & Z^{'} \ar[d]_{w'} \\ Z \ar[r]_{v}& Z \ar[r]_{w} & Y}. \hspace*{1cm}(*)$$
Indeed, by the statement of the proposition, the left square is homotopy cartesian
and now we can form the homotopy pushout of
$\xymatrix{Z&\ar[l]_{t}N\ar[r]^-{u's}& Z'}$,
completing with morphisms
$\xymatrix{Z\ar[r]^{\eta}&W&\ar[l]_{\rho}Z'}$ such that we get the following homotopy cartesian square
$$
\xymatrix{
N\ar[d]_t\ar[r]^{u's}&Z'\ar[d]_{\rho}\\
Z\ar[r]_\eta&W
.}
$$
Since we have an equality $\eta t=\rho u' s$, by properties of homotopy
pushouts, we get a not necessarily unique morphism $\varphi :Z\longrightarrow W$
such that $\varphi v=\eta$ and $\varphi u=\rho u'$. This leads to the following
commutative diagram, where the left square and the patching of the two
squares are homotopy cartesian:

$$\xymatrix{N \ar[d]_{t} \ar[r]^{s} & M \ar[d]_{u} \ar[r]^{u'} & Z^{'} \ar[d]_{\rho} \\ Z \ar[r]_{v}& Z \ar[r]_{\varphi} & W}. $$
By \cite[Lemma 9]{SZ2}, one can replace $\varphi$ by an appropriate substitute so
that the right square is also homotopy cartesian.  But then $W$ is isomorphic to
$Y$ and, looking at this isomorphism as an identification, we can choose $w'=\rho$
and $w$ to be the mentioned substitute of $\varphi$.

Now the different homotopy cartesian squares fit as faces of the following almost-cube:
$$
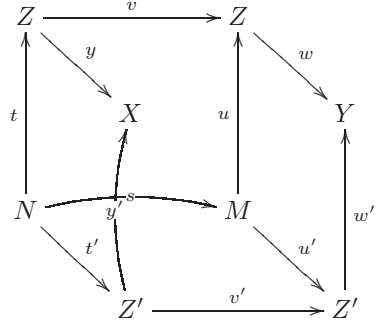
\xymatrix{Z\ar[rr]^v\ar[dr]^y&&Z\ar[rd]^w\\
&X&&Y\\
N\ar[uu]^t\ar@/^/[rr]|-s\ar[dr]^{t'}&&M\ar[uu]^u\ar[dr]^{u'}\\
&Z'\ar[rr]^{v'}\ar@/^/[uu]|-{y'}&&Z'\ar[uu]_{w'}
}
$$

Since $X$ is a homotopy pushout, and since
$$wvt=wus=w'u's=w'v't',$$
there is a morphism $X\stackrel{d}{\lra} Y$ such that
$$yd=wv\text{ and }dy'=w'v'.$$
We can complete the vertical morphisms to distinguished triangles and obtain hence the
commutative diagram with vertical sequences being distinguished triangles.
$$
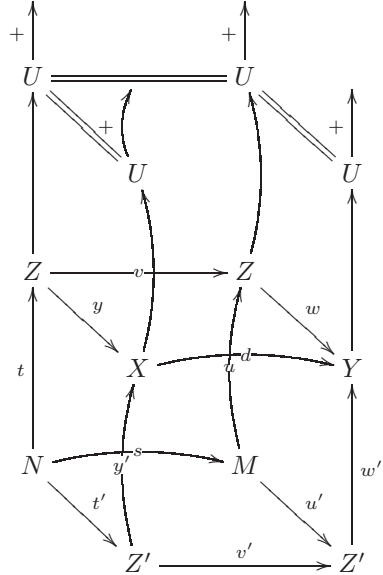
\xymatrix{&&&\\
U\ar@{=}[rr]\ar@{=}[dr]\ar[u]^+&&U\ar@{=}[dr]\ar[u]^+&&\\
&U\ar@/^/[u]^+&&U\ar[u]^+\\
Z\ar[rr]|-v\ar[dr]^y\ar[uu]&&Z\ar[rd]^w\ar@/_/[uu]\\
&X\ar@/^/[rr]|-d\ar@/_/[uu]&&Y\ar[uu]\\
N\ar[uu]^t\ar@/^/[rr]|-s\ar[dr]^{t'}&&M\ar@/^/[uu]|-u\ar[dr]^{u'}\\
&Z'\ar[rr]^{v'}\ar@/^/[uu]|-{y'}&&Z'\ar[uu]_{w'}
}
$$

In particular, we get the following commutative diagram:

$$\xymatrix{N \ar[d]_{t} \ar[r]^{t'} & Z^{'} \ar[d]_{y'} \ar[r]^{v'} & Z^{'} \ar[d]_{w'} \\ Z \ar[r]_{y}& X \ar[r]_{d} & Y}.$$
Its left square is  homotopy cartesian by hypothesis and the patching
of the two square is the pathching of the two squares in the diagram (*)
above, and then is also homotopy cartesian.
 By \cite[Lemma 9]{SZ2} again,  replacing $d$ by a suitable $d'$, we can
 assume without loss of generality that the right square is also homotopy cartesian.
But this shows that $X\leq_{\Delta,\textup{left}}Y$ (resp.
$X\leq_{\Delta+\text{nil,left}}Y$ when $v$ and $v'$ are nilpotent).
\end{proof}

\begin{Cor}
Suppose that
$\mathcal T$ has split idempotents and is either the category of compact objects
in an algebraic compactly generated triangulated category, or else has artinian
endomorphism rings of objects.
In the notation of Proposition~\ref{pushouts}, then
$$M\leq_{\Delta+\textup{nil}}N\Rightarrow X\leq_{\Delta+\textup{nil}}Y.$$
\end{Cor}

\begin{proof}
Indeed, using \cite{SZ2}, we see that
$X\leq_{\Delta+\textup{nil,right}}Y\Leftrightarrow X\leq_{\Delta+\text{nil,left}}Y$
under these hypotheses. \end{proof}

\section{Degeneration of zero and the zero objects in the Grothendieck group}

\label{ZeroinGrothendieck}

Recall that, given full subcategories $\mathcal{U}$ and $\mathcal{V}$ of a triangulated
category $\mathcal{T}$, then the subcategory $\mathcal{U}\star\mathcal{V}$ is the full
subcategory of $\mathcal T$ consisting of the objects $M$ that fit in a
distinguished triangle
$U\longrightarrow M\longrightarrow V\stackrel{}{\longrightarrow}U[1]$, with $U\in\mathcal{U}$
and $V\in\mathcal{V}$. It is well-known that the operation $\star$ is associative, in the
sense that
$(\mathcal{U}\star\mathcal{V})\star\mathcal{W}=\mathcal{U}\star (\mathcal{V}\star\mathcal{W})$,
for all subcategories $\mathcal{U}, \mathcal{V},\mathcal{W}$ of $\mathcal{T}$
(see \cite[Lemme 1.3.10]{BBD}). If one puts
$\mathcal{U}^{\star n}=\underbrace{\mathcal{U}\star{\cdots}\star\mathcal{U}}_{n\text{ factors}}$,
for each $n\geq 0$ (with the convention that $\mathcal{U}^{\star 0}=0$), then
$\mathcal{U}^{\text{ext}}=\bigcup_{n\in\mathbf{N}}\mathcal{U}^{\star n}$ is the
\emph{extension closure} of $\mathcal{U}$, that is, the smallest
subcategory of $\mathcal{T}$
closed under extensions that contains $\mathcal{U}$. The smallest
triangulated subcategory of
$\mathcal{T}$ that contains $\mathcal{U}$, denoted
$\text{tria}_\mathcal{T}(\mathcal{U})$, is
$$\text{tria}_\mathcal{T}(\mathcal{U})=
\bigcup_{n\in\mathbb{N}}\bigcup_{(r_1,...,r_n)\in\mathbb{Z}^n}\mathcal{U}[r_1]\star \cdots\star\mathcal{U}[r_n].$$
In other words, the objects of  $\text{tria}_\mathcal{T}(\mathcal{U})$ are
precisely those $M$ admitting a sequence $$0=M_0\stackrel{f_1}{\longrightarrow}M_1\stackrel{f_2}{\longrightarrow}\cdots \stackrel{f_{n-1}}\lra M_{n-1}\stackrel{f_n}{\longrightarrow}M_n=M,$$
where $\text{cone}(f_k)$ is isomorphic to $U_k[r_k]$, for some $U_k\in\mathcal{U}$ and some $r_k\in\mathbb{Z}$, for all $k=1,...,n$.

Let $\mathcal T$ be a skeletally small triangulated category.  All throughout
this section we will put $\leq_\Delta =\leq_{\Delta, \textup{left}}$ and
$\leq_{\Delta+\textup{nil}} =\leq_{\Delta+\textup{nil ,left}}$
Recall from Proposition~\ref{JSZ2prop} that $\leq_{\Delta}$ can be shown to
be a transitive relation in some cases. On the other hand,
 we will denote by  $\preceq_{\Delta+\textup{nil}}$ the
smallest transitive relation containing $\leq_{\Delta+\textup{nil}}$.

\begin{Theorem} \label{thm.zero objects in K0}
Let $\mathcal{S}$ be a set of objects in the triangulated category $\mathcal{T}$
such that $\mathcal{T}=\textup{tria}_\mathcal{T}(\mathcal{S})$, let
$[\mathcal{S}]:=\{[S]\text{: }S\in\mathcal{S}\}$ denote the corresponding
set of generators of the group $K_0(\mathcal{T})$ and let $\widehat{\mathcal{S}}$
be the subcategory of $\mathcal{T}$ consisting of the objects $X$ which
are finite direct sums of shifts of  objects in $\mathcal{S}$ and are such that $[X]=0$ in
$K_0(\mathcal{T})$. Denote by
\begin{itemize}
\item
$\mathcal{T}^0_\Delta$
(resp. $\mathcal{T}^0_{\Delta +\textup{nil}}$) the full subcategory of
$\mathcal{T}$ consisting of the objects $X$ such that
$0\leq_\Delta X$ (resp. $0\leq_{\Delta +\textup{nil}}X$)
\item and by
$\mathcal{T}^0$ the (triangulated) subcategory of $\mathcal{T}$
consisting of the objects $M$ such that $[M]=0$ in the group $K_0(\mathcal{T})$.
\end{itemize}
Then the following assertions hold:
\begin{enumerate}
\item An object $M$ is in $\mathcal{T}^0$ if, and only if,
$M\leq_{\Delta}X$ (resp. $M\preceq_{\Delta+\textup{nil}}X$),  for some $X\in\widehat{\mathcal{S}}$.
When $[\mathcal{S}]$ is a basis of $K_0(\mathcal{T})$ the objects $X$
can be chosen to be finite direct sums of shifts of objects in
$\bar{\mathcal{S}}:=\{S\oplus S[2k+1]\text{: }k\in\mathbb{Z}; S\in{\mathcal S}\}$.
\item $\mathcal{T}^0=\text{tria}_\mathcal{T}(S\oplus S[t_S]\text{: }S\in\mathcal{S})$,
for  every choice of odd integers $t_S$.
\item $\mathcal{T}^0$ is the extension closure of $\mathcal{T}_\Delta^0$
(resp. $\mathcal{T}_{\Delta +\textup{nil}}^0$).
\end{enumerate}
\end{Theorem}

\begin{proof}
(1) By Corollary \ref{K0triangledeg}, the `if' part of this implication is clear.
For the `only if' part, we first claim that, for each $M\in\mathcal{T}$, one has
that
$M\preceq_{\Delta+\textup{nil}}\oplus_{S\in \mathcal{S}}\oplus_{k\in\mathbb{Z}}S[k]^{m_{S,k}}$,
where the $S$ are in $\mathcal{S}$ and the $m_{k,S}$ are nonnegative integers,
all zero but a finite
number. Recall that $\mathcal{T}=\text{tria}_\mathcal{T}(\mathcal{S})$,
and so we have a finite sequence
$$0=M_0\stackrel{f_1}{\longrightarrow}M_1\stackrel{f_2}{\longrightarrow}\dots
\stackrel{f_n}{\longrightarrow}M_n=M\hspace*{1cm} (*)$$
such that $C_k:=\text{cone}(f_k)$ is a shift of some object
of $\mathcal{S}$, for each $k=1,...,n$. We will settle our
claim by induction on $n>0$,  the case $n=1$ being clear.
Suppose now that $n>0$ and consider the induced triangle
$$\xymatrix{M_{n-1}\ar[r]^-{f_n}&M\ar[r]^-{g}&C_n
\ar[r]&M_{n-1}[1],}$$
where $C_n\cong S[k]$, for some $S\in \mathcal{S}$ and $k\in\mathbb{Z}$.
Taking the homotopy pushout of $f_n$ and the zero endomorphism
$M_{n-1}\stackrel{0}{\longrightarrow}M_{n-1}$, we readily see that we
have a distinguished triangle
$$\xymatrix{M_{n-1}\ar[r]^-{\scriptsize\begin{pmatrix} 0\\ f_n\end{pmatrix}}&M_{n-1}\oplus M\ar[r]& M_{n-1}\oplus C_n\ar[r]&M_{n-1}[1].}$$
That is, we have $M\preceq_{\Delta +\textup{nil}}M_{n-1}\oplus C_n\cong M_{n-1}\oplus S[k]$.
The result then follows by induction since $A_i\preceq_{\Delta+\textup{nil}} B_i$, for $i=1,2$,
implies that $A_1\oplus A_2\preceq_{\Delta+\textup{nil}} B_1\oplus B_2$.

We also claim that
$M\leq_\Delta\oplus_{S\in \mathcal{S}}\oplus_{k\in\mathbb{Z}}S[k]^{m_{S,k}}$,
for $S_k$ and $m_{S,k}$
as in the previous paragraph.  Using again the sequence $(*)$ and bearing in
mind that each cone $C_k$ is a shift of some object in $\mathcal{S}$, we
consider the distinguished triangles
$$\xymatrix{M_{k-1} \ar[r]^-{f_k}&M_k\ar[r]& C_k\ar[r]&M_{k-1}[1]}$$
for all $k\in\{1,\dots,n-1\}$. Taking the direct sum of these
distinguished triangles we get a distinguished triangle
$$\xymatrix{\left(\bigoplus_{k=1}^{n-1}M_k\right)\ar[r]^-{\oplus_{k=1}^nf_k}& \left(M\oplus \bigoplus_{k=1}^{n-1}M_k\right)\ar[r]& \left(\bigoplus_{k=1}^{n}C_k\right)\ar[r]&\left(\bigoplus_{k=1}^{n-1}M_k\right)[1]}$$
and hence $M\leq_\Delta\bigoplus_{k=1}^nC_k$, as desired.

The last two paragraphs show that we have $M\preceq_{\Delta +\textup{nil}}X$
and $M\leq_\Delta Y$, for objects $X,Y$ which are direct sums of shift of
objects of $\mathcal{S}$.
When in addition $M\in\mathcal{T}^0$, by Corollary  \ref{K0triangledeg}, we also
have $[X]=[Y]=0$ in $K_0(\mathcal{T})$. Therefore we have that
$X,Y\in\hat{\mathcal{S}}$. This proves assertion 1, except for the final statement.

To prove that final statement, suppose that $[\mathcal{S}]$ is a
basis of $K_0(\mathcal{T})$. We shall prove that in this case each
object of $\hat{\mathcal{S}}$ is a direct sum  of objects of the
form $S[k]\oplus S[l]=(S\oplus S[l-k])[k]$, with $S\in\mathcal{S}$
and exactly one of $l$ and $m$ being odd. This will end the proof.
Let then take $X\in\hat{\mathcal{S}}$ and decompose it as
$X=\oplus_{S\in \mathcal{S}}\oplus_{k\in\mathbb{Z}}S[k]^{m_{S,k}}$.
 Note that,  due to the fact that  $[\mathcal{S}]$ is a basis of
 $K_0(\mathcal{T})$,
the summand $X_S=\oplus_{k\in\mathbb{Z}}S[k]^{m_{k,S}}$ also
satisfies that $[X_S]=0$
in $K_0(\mathcal{T})$, for each $S\in\mathcal{S}$. So it is not restrictive to
assume that $X=S[k_1]^{m_{1}}\oplus S[k_2]^{m_2}\oplus \cdots\oplus S[k_r]^{m_r}$,
for some pairwise different  integers $k_1,\dots,k_r$, where, for simplicity, we
have put $m_{k_i,S}=m_i>0$ for $i=1,\dots,r$.  We can reorder the summands in this
last direct sum, so that $k_i$ is even, for $1\leq i\leq t$, and $k_i$ is odd,
for $t<i\leq n$. Bearing in mind that $[S[k]]=(-1)^k[S]$ in $K_0(\mathcal{T})$,
that $[\mathcal{S}]$ is a basis of $K_0(\mathcal{T})$ and that $[X]=0$ in this
latter abelian group, we conclude that $\sum_{i=1}^tm_i=\sum_{i=t+1}^nm_i$.
We call $m(X)$ this last integer which is strictly positive when $X\neq 0$.

We now prove the result by induction on $m(X)>0$. If $m(X)=1$ then we have
that $X\cong S[k]\oplus S[l]$, where $k$ is even and $l$ is odd, and we are done.  If $m>1$,
we put $q=\text{min}\{m_1,\dots,m_r\}$. We decompose $m_1=q+m'_1$ and $m_n=q+m'_n$
and $m'_i=m_i$, for all $i\neq 1,n$. Then we have a decomposition
$X=(S[k_1]\oplus S[k_n])^q\oplus X'$, where $k_1$ and $k_n$ are even and
odd, respectively, and $X'= \oplus_{i=1}^nS[k_i]^{m'_i}$ is either zero or a nonzero summand such that $[X']=0$ in $K_0(\mathcal{T})$ and $m(X')<m(X)$. Then the induction hypothesis applies.

(2) Let $(t_S)_{S\in\mathcal{S}}$ be a collection of odd integers and put $\mathcal{D}:=\text{tria}_\mathcal{T}(S\oplus S[t_S]\text{: }S\in\mathcal{T}_S)$.
It follows that each object of $\mathcal{S}$ is a direct summand of an object of
$\mathcal{D}$ and since each object $T$ of
$\mathcal{T}=\text{tria}_\mathcal{T}(\mathcal{S})$ is a finite iterated
extension of objects of the form $S[k]$, with $S\in\mathcal{S}$ and
$k\in\mathbb{Z}$, it easily follows that each such $T$ is a direct summand
of an object of $\mathcal{D}$. This means that $\mathcal{D}$ is a dense
triangulated subcategory of $\mathcal{T}$ in the terminology of \cite{Thomason}.
Moreover, we clearly have $\mathcal{D}\subseteq\mathcal{T}^0$. But
\cite[Theorem 2.1]{Thomason} gives an order-preserving bijection between
the dense triangulated subcategories of $\mathcal{T}$ and the subgroups of
$K_0(\mathcal{T})$. Since $\mathcal{T}^0$ corresponds to $0$ by this
bijection we get that $\mathcal{D}=\mathcal{T}^0$, as desired.

(3)
Note that assertion 2 implies assertion 3.
Indeed, by the comments preceding Theorem~\ref{thm.zero objects in K0},
assertion 2 says that $\mathcal{T}^0$ is the extension closure of
$\{(S\oplus S[t_S])[n]\text{: }S\in\mathcal{S}\text{ and }n\in\mathbb{Z}\}$,
for any choice of odd integers $t_S$ ($S\in\mathcal{S}$). We may choose $t_S=1$
for each $S$, and then it is obvious that $0\leq_{\Delta+\text{ nil}}S\oplus S[1]$.
Since we have inclusions
$$\{(S\oplus S[1])\text{: }S\in\mathcal{S}\}\subset\mathcal{T}_{\Delta+\textup{nil}}^0
\subset\mathcal{T}_\Delta^0\subseteq\mathcal{T}^0$$
assertion 3 immediately follows.
\end{proof}

\begin{Example}
The following are examples of a triangulated category $\mathcal{T}$ and a set
$\mathcal{S}$ of its objects that satisfy the hypotheses of
Theorem~\ref{thm.zero objects in K0} and, in addition, $[\mathcal{S}]$
is a basis of $K_0(\mathcal{T})$. Here
$K$ is a commutative ring

\begin{enumerate}
\item Call a dg $K$-algebra $A$ \emph{homologically non positive} when
$H^kA=0$, for all $k>0$, and call it \emph{homologically finite dimensional}
when $H^*(A)=\oplus_{k\in\mathbb{Z}}H^k(A)$ is a $K$-module of finite length.
For instance, any Artin algebra is homologically non positive and homologically
finite dimensional over its center, when viewed as dg algebra. Let $A$ be a
homologically non positive homologically finite dimensional dg algebra and let
$\mathcal{T}=\mathcal{D}^b_{fl}(A)$ be the subcategory of the derived category
$\mathcal{D}(A)$ consisting of the dg $A$-modules $M$ such that
$H^*(M)=\oplus_{k\in\mathbb{Z}}H^k(M)$ has finite length as a $K$-module.
When choosing as $\mathcal{S}$ a set of representatives, up to isomorphism
in $\mathcal{D}(A)$, of the dg $A$-modules $S$ such that $H^*(S)=H^0(S)$
(i.e. its homology is concentrated in degree zero) and $H^0(S)$ is a
simple $H^0(A)$-module, one has that $\mathcal{T}$ and $\mathcal{S}$
satisfy the hypotheses of  Theorem~\ref{thm.zero objects in K0}
and $[\mathcal{S}]$ is a basis of $K_0(\mathcal{T})$. In particular, taking $A$
to be an Artin algebra, $\mathcal{T}=\mathcal{D}^b(\textup{mod}-A)$ and
$\mathcal{S}$ be a set of representatives, up to isomorphism, of the
simple $A$-modules (viewed as stalk complexes in degree zero), the
hypotheses of  Theorem~\ref{thm.zero objects in K0}  hold and
$[\mathcal{S}]$ is a basis of $K_0(\mathcal{T})$.

\item Suppose that $\mathcal{A}$ is an additive category
with a set $\mathcal{S}'$ of objects such that
$\mathcal{A}=\text{add}(\mathcal{S}')$ and the Grothendieck group
$K_0(\mathcal{A})$ is free with $\{[S]\text{: }S\in\mathcal{S}'\}$ as a basis.
Then the bounded homotopy category
$\mathcal{T}=\mathcal{K}^b(\mathcal{A})$ and the set
$\mathcal{S}=\mathcal{S}'[0]$ of stalk complexes
$S'[0]$, with $S'\in\mathcal{S}'$, satisfy the hypotheses
of Theorem~\ref{thm.zero objects in K0}
(see \cite[Theorems 1.1 and 1.2]{Rose}). This includes the case when
$\mathcal{T}=\mathcal{K}^b(A-\text{proj})$, where $A$ is
a principal ideal domain or a semiperfect ring, in particular an Artin algebra,
by taking as $\mathcal{S}'$ the set of (isomorphism classes of)
indecomposable projective $A$-modules
\end{enumerate}
\end{Example}

\begin{Lemma}  \label{counterexampleforDelta}
Let $A$ be an Artin algebra and $S$ be a simple $A$-module.
For each integer $k\neq -1,0$, the complex $M=S\oplus S[2k+1]$ has the
property that $[M]=0$ in $K_0(D^b(A-\textup{mod}))$, but it is not a
$\Delta$-degeneration of zero:
$0\not\leq_\Delta M$.
\end{Lemma}

\begin{proof}
Note that the homology module $H^i(M)$ is zero,
except for $i=0$ and $i=-2k-1$. If there is a distinguished triangle
$$\xymatrix{Z\ar[r]^{f}&Z\ar[r]& M\ar[r]&Z[1]}$$
in $D^b(A-\text{mod})$, the associated sequence of homologies
gives an exact sequence
$$\xymatrix{0\ar[r]&H^0(Z)\ar[r]^-{H^0(f)}&H^0(Z)\ar[r]& S\ar[r]& H^1(Z)\ar[r]^-{H^1(f)}&H^1(Z)\ar[r]& 0,}$$
which forces $H^0(f)$ and $H^1(f)$ to be isomorphisms since
they are a monomorphic and an epimorphic endomorphism,
respectively, of finite length modules. Therefore $S=0$, and we get a contradiction.
\end{proof}

\begin{Rem} Using Theorem~\ref{SZ1thm} and Proposition~\ref{JSZ2prop} we see that
in the situation of Lemma~\ref{counterexampleforDelta}
the relations ``$\leq_{\Delta+\textup{ nil}}$'', ``$\leq_\Delta$'',
and ``$\leq_{\textup{cdeg}}$'' coincide.
\end{Rem}

Recall that,  for a skeletally small
triangulated category $\mathcal T$,  we denote by $\preceq_{\Delta+\textup{ nil}}$ the
smallest transitive relation on the set of isomorphism classes of objects in $\mathcal T$
containing $\leq_{\Delta+\textup{ nil}}$.

\begin{Prop}
\label{Counterexample2implies1-b}
Let  $\mathcal{A}$ be any skeletally small abelian category for which
$\mathcal{D}^b(\mathcal{A})$ is well-defined, i.e. it has $\text{Hom}$
sets as opposed to proper classes, and let us identify $\mathcal{A}$ with the subcategory of  $\mathcal{D}^b(\mathcal{A})$ consisting of objects $X$ such that $H^i(X)=0$,
for $i\neq 0$. The following assertions hold:
\begin{enumerate}
\item If $Y\leq_{\Delta+\textup{nil}}X$ in $\mathcal{D}^b(\mathcal{A})$ and $X\in\mathcal{A}$,
then $Y\in\mathcal{A}$.
\item If there is an object in $\mathcal{A}$
which has a monomorphic endomorphism which is not an isomorphism, then there
is an object $X$
in $\mathcal{T}:=\mathcal{D}^b(\mathcal{A})$  such that $0\leq_\Delta X$
(and hence $[X]=0$ in $K_0(\mathcal{T})$), but $0\not\preceq_{\Delta+\textup{nil}}X$.
\end{enumerate}

Assertion 2 includes in particular the case of the category $\mathcal{A}=R-\text{mod}$
of finitely generated $R$-modules, when  $K$ is a
integral domain and $R$ is a noetherian $K$-algebra which is torsion-free as a $K$-module
and does not contain the field of fractions of $K$ as a subalgebra.
\end{Prop}

\begin{proof}
(1) Let us consider a distinguished triangle
$$\xymatrix{W\ar[r]^-{\scriptsize\begin{pmatrix}v\\ \alpha \end{pmatrix}}&W\oplus Y\ar[r]& X\ar[r]& W[1]}$$
 in $D^b(\mathcal{A})$, where $v$ is a
nilpotent endomorphism of $W$ and $X\in\mathcal{A}$.  The long exact sequence of homologies
gives that
$$\xymatrix{H^j(W)\ar[rr]^-{\scriptsize\begin{pmatrix} H^j(v)\\ H^j(\alpha) \end{pmatrix}}&&H^j(W)\oplus H^j(Y)}$$
is an isomorphism, for $j\neq 0,1$, and there is an exact sequence
$$\xymatrix{0\ar[r]& H^0(W)\ar[r]^-{\scriptsize\begin{pmatrix} H^0(v)\\ H^0(\alpha) \end{pmatrix}}&
H^0(W)\oplus H^0(Y)\ar[r]& X\ar[r]& H^{1}(W)\ar[r]^-{\scriptsize\begin{pmatrix} H^1(v)\\ H^1(\alpha) \end{pmatrix}}&H^1(W)\oplus H^1(Y)\ar[r]& 0}$$
in $\mathcal{A}$.
Proving that $Y$ has homology concentrated in zero degree reduces to prove that if
$\begin{pmatrix}w\\ g \end{pmatrix}:A\longrightarrow A\oplus B$ is an
epimorphism in $\mathcal{A}$, for some objects $A,B\in\mathcal{A}$,
where $w$ is a nilpotent endomorphism of $A$, then $A=B=0$.
This is clear when $w=0$.
But if $w\neq 0$ and $m$ is the nilpotent index of $w$
(i.e. $w^m=0\neq w^{m-1}$), then the composition
$$\xymatrix{A\ar[r]^{\scriptsize\begin{pmatrix}w\\ g \end{pmatrix}}&A\oplus B\ar[rr]^{\scriptsize\begin{pmatrix}w^{m-1}& 0 \end{pmatrix}}&&A}$$
is the zero map, which implies that $w^{m-1}=0$,
thus yielding a contradiction.

\bigskip

(2) Let $f:Z\longrightarrow Z$ be a monomorphic endomorphism
which is not an isomorphism and put $X=\textup{coker}(f)$.
We then have an induced  distinguished triangle
$$\xymatrix{Z\ar[r]^-{f}&Z\ar[r]& X\ar[r]& Z[1]}$$
in $\mathcal{D}^b(\mathcal{A})$, thus showing that $0\leq_\Delta X=X[0]$
in the later triangulated category.
Suppose now that  $0\preceq_{\Delta+\textup{nil}}X$. Then
we have a sequence $0=X_0,X_1,\dots,X_n=X$ in $D^b(\mathcal{A})$ such
that $X_{i-1}\leq_{\Delta+\textup{nil}} X_i$ and $X_i\neq 0$ for $i=1,\dots,n$.
By assertion 1, we know that all $X_i$ are in $\mathcal{A}$.
Replacing $X$ by $X_1$ if necessary, we get an object $X\neq 0$ of $\mathcal{A}$
such that $0\leq_{\Delta+\textup{nil}}X$ in $\mathcal{D}^b(\mathcal{A})$.
We can fix a distinguished triangle
$$\xymatrix{Q\ar[r]^-{u}&Q\ar[r]&X\ar[r]&Q[1]}$$ in
$\mathcal{D}^b(\mathcal{A})$, where $u$ is a nilpotent endomorphism of $Q$.
The long exact sequence of homologies gives then an exact sequence
$$\xymatrix{0\ar[r]& H^0(Q)\ar[r]^-{H^0(u)}&H^0(Q)\ar[r]& X\ar[r]&
H^1(Q)\ar[r]^-{H^1(u)}&H^1(Q)\ar[r]&0}$$
in $\mathcal{A}$. But it is obvious that a nilpotent endomorphism of an object
$A'\in\mathcal{A}$ can be a monomorphism or an epimorphism only in case $A'=0$.
We then get $H^j(Q)=0$ for $j=0,1$, which in turn implies $X=0$ and hence a contradiction.

Finally, suppose now that $R-\text{mod}$ is the category of finitely generated
modules over a noetherian $K$-algebra $R$ as indicated in the last statement of
the proposition. Since the field of fractions of $K$ is not a subring of $R$,
there must be a $\lambda\in K\setminus\{0\}$ such that $\lambda R\subsetneq R$.
Due to the torsion-free condition of $R$ as a $K$-module,  multiplication by
$\lambda$ gives a monomorphism
$u=u_\lambda :R\longrightarrow R$ in $R-\textup{mod}$ which is not an
epimorphism.
\end{proof}

\begin{Example}\label{Counterexample2implies1}
Note that Proposition~\ref{Counterexample2implies1-b} shows that
the stalk complex $\Z/p\Z$ in $\mathcal{D}^b(\Z-\textup{mod})$ has zero image in
the Grothendieck group of the derived category and that $0\not\preceq_{\Delta+\textup{ nil}}\Z/p\Z$.
\end{Example}

\begin{Rem}
Proposition~\ref{Counterexample2implies1-b} shows that
being degeneration of zero is strictly stronger than having zero image
in the Grothendieck group.
\end{Rem}

\end{document}